\documentclass[12pt,reqno]{article}
\usepackage{amsmath}
\usepackage{amssymb}
\usepackage{amsthm}

\usepackage{algorithmic}

\usepackage{pstricks}
\usepackage{multido}

\usepackage{graphicx}
\usepackage{color}

\newcommand{\showgrid}{}

\newcommand{\gridon}{\renewcommand{\showgrid}{\psset{subgriddiv=1,griddots=10,gridlabels=6pt}\psgrid}}

\gridon 

\def\figref#1{Figure~\ref{#1}}
\def\heute{\number\day.~\ifcase\month\or
 J\"anner\or Februar\or M\"arz\or April\or Mai\or Juni\or
 Juli\or August\or September\or Oktober\or November\or Dezember\fi
 \space\number\year}

\def\defeq{:=}
\def\absof#1{\left|#1\right|}
\def\setof#1{\left\{#1\right\}}
\def\of#1{\!\left(#1\right)}

\def\pas#1{\left(#1\right)}
\def\brk#1{\left[#1\right]}
\def\sgn{{\operatorname{sgn}}}

\def\bit{\begin{itemize}}
\def\eit{\end{itemize}}
\def\beq{\begin{equation}}
\def\eeq{\end{equation}}

\def\R{{\mathbb R}}
\def\Z{{\mathbb Z}}

\def\1{{\mathbf 1}}

{\begin{list}{#1}{\parsep0em \itemsep0.3em \labelwidth1em
\labelsep0.5em \leftmargin1.5em }} {\end{list}}

\def\path{{%\mathrm 
P}}
\def\figref#1{Figure~\ref{#1}}

\def\EM#1{{\em #1\/}}

\newif\ifenglish
\englishtrue

\newtheorem{thm}{\ifenglish Theorem\else Satz\fi}

\theoremstyle{remark}

\def\pfaffian{\operatorname{Pf}}

\def\weight{\omega}

\def\and{\wedge}

\def\myrange#1{\setof{1,\dots,k}}

\def\strich{^\prime}

\def\ueberschrift{\subsubsection*}

\newcount\labelcount

\newif\iflongversion
\longversionfalse

\newif\ifverylongversion
\verylongversionfalse

\begin{document}

\bibliographystyle{plain}

\title{Speyer's elegant topological proof for Kasteleyn's Theorem}

\begin{abstract}
The purpose of this note is to rephrase Speyer's elegant topological
proof for Kasteleyn's Theorem in a simple graph theoretical manner.
\end{abstract}

\author{Markus Fulmek}

\thanks{
Research supported by the National Research Network ``Analytic
Combinatorics and Probabilistic Number Theory'', funded by the
Austrian Science Foundation. 
}

\maketitle

\psset{linewidth=0.05,linecolor=gray,linearc=0.3} 
\def\drawvertex#1{\pscircle[linewidth=0.05,linecolor=black,fillstyle=solid,fillcolor=black](#1){0.2}}
\def\drawnonvertex#1{\pscircle[linewidth=0.05,linecolor=black,fillstyle=solid,fillcolor=white](#1){0.175}}

Speyer \cite{Speyer:2016:VOATOKW} recently published a very short and elegant proof
for Kasteleyn's Theorem \cite{Kasteleyn:1967:GTACP}. This proof involves some higher topology arguments; and
the purpose of this note is to replace them by elementary graph--theoretical
arguments. 
We shall, however, present only a sketch of a proof, since the technical details
of a fully rigorous proof would likely conceal the beauty and simplicity of the basic idea.


\section{Basic definitions and notation}

For reader's convenience and to fix our notation, we shall briefly
recall some basic definitions.

\ueberschrift{Graph}

A finite graph $G$ consists
\bit
\item of a finite set of
\EM{vertices} $V\of G$
\item and a finite (multi--)set of \EM{edges} $E\of G$,
\eit
where each $e\in E\of G$
\EM{joins} two different vertices $v_1, v_2\in V\of G$ (we also say that $e$ is \EM{incident}
with $v_1$ and $v_2$), which we denote by $e=\setof{v_1,v_2}$.
We shall not consider \EM{loops}, which are
edges ``joining a single vertex'' ($e=\setof{v}$): If such loop arises during the constructions
described below, we shall remove it immediately.

\ueberschrift{Graph embedding}

A finite graph $G$ can always be \EM{drawn} in the plane $\R^2$, such that
\bit 
\item the vertices correspond bijectively to \EM{points} in the plane (we shall call them \EM{vertex--points}),
\item and the edges correspond bijectively to smooth \EM{curves} $\brk{0,1}\to\R^2$ in the plane (we shall call them \EM{edge--curves})
	which \EM{connect} the respective vertex--point, i.e.,
	some edge $e=\setof{v_1,v_2}$ corresponds to an edge--curve
	\bit
	\item which starts in the vertex--point corresponding to $v_1$
	and ends in the vertex--point corresponding to $v_2$ (or vice versa: the orientation of the curve is irrelevant in our context),
	\item but \EM{does not} ``intersect itself'' (i.e., does not pass twice through the same point;
	in particular, edge--curves do not ``return'' to their starting point or end point),
	\item and does not pass through any other vertex--point.
	\eit
\eit
Such drawing is called an \EM{embedding} of the graph $G$ in the plane.

By definition, the edge--curves of edges incident to the same vertex $v$ have
the vertex--point corresponding to $v$ in common,
but there might be other such intersections in \EM{nonvertex--points}.

For our purpose, we additionally assume that:
\begin{enumerate}
\item Any two curves intersecting in some nonvertex--point $p$
	do actually \EM{cross} each other (i.e., they are \EM{not tangent} in $p$),
	whence we call such point $p$ a \EM{crossing} of the two curves.
\item There is only a finite number of such \EM{crossings}.
\end{enumerate}
We call an embedding which fulfils these additional assumptions a \EM{proper embedding}.

See \figref{fig:simple-graph} for an
illustration of these concepts: In all graphical representations of proper embeddings in this
note, we shall indicate
\bit
\item vertex--points by \EM{black dots},
\item curves corresponding to edges by \EM{gray lines},
\item and crossings by \EM{small white circles}.
\eit
\begin{figure}
\label{fig:simple-graph}
\caption{Example of a proper embedding of some graph $G$ in the plane.}
{\small
Let $V\of G = \setof{v_1,v_2,\dots, v_6}$ and $E\of G = \setof{
\setof{v_1,v_2},
\setof{v_2,v_3}^2,
\setof{v_3,v_4},
\setof{v_5,v_6}}$, where the notation $\setof{v_2,v_3}^2$ indicates that there are
\EM{two} edges joining vertices $v_2$ and $v_3$. The picture shows a proper embedding of $G$
with a single crossing $p$: The points corresponding to the vertices are indicated by black
dots, the curves corresponding to edges are indicated by gray lines, and the crossing is indicated by a small
white circle.
}
\begin{center}
\input graphics/simple-graph.tex
\end{center}
\end{figure}

\ueberschrift{Simplified notation}
If we are given some graph $G$ and an embedding $\eta$ of $G$, let us denote
\bit
\item by $\eta\of v$ the vertex--point representing vertex $v\in V\of G$,
\item by $\eta\of e$ the edge--curve representing edge $e\in E\of G$.
\eit
But we shall abuse this notation in the following by making no distinction
\bit
\item between $v\in V\of G$ and $\eta\of v$, i.e., instead of saying
	``the vertex--point in the plane representing $v$'' we simply shall say ``$v$'',
\item and between $e\in E\of G$ and $\eta\of e$, i.e., instead of saying
	``the edge--curve in the plane representing $e$'' we simply shall say ``$e$''.
\eit
This slight imprecision should cause no confusion in our context.

\ueberschrift{Planar graph}

A proper embedding without any intersection of edge--curves in non--vertex points is called a
\EM{planar embedding}.
A graph which admits a planar embedding is called a \EM{planar graph}.

Not all graphs are planar, but we may \EM{view} every proper embedding of some 
graph $G$ as planar embedding of a planar graph $G\strich$ by \EM{reinterpreting} all
crossings as vertex--points: $G=G\strich$ (in the sense of graph isomorphisms) if the
embedding has no crossings at all. (In \figref{fig:simple-graph}, this reinterpretation would amount
to replacing the single crossing $p$ by a new vertex--point $v_7$.)

\ueberschrift{Perfect matching}

A \EM{perfect matching} $M\subseteq E\of G$ of $G$ is a subset of edges
such that every vertex of $G$ is incident with precisely one edge in $M$.

\ueberschrift{Edge weight}

We assume that we are given some (nowhere--zero) \EM{weight function} (called \EM{edge weight}) $\weight$
$$\weight:E\of G\to R\setminus\setof{0}$$
on the multiset of edges of $G$, where $R$ is some (nontrivial) commutative
ring (in most cases, $R$ is $\Z$ or some polynomial ring; the constant edge weight 
$\weight\equiv 1$ is used for enumeration purposes).

%

\ueberschrift{Generating function of perfect matchings}

The weight of some perfect matching $M\subseteq E\of G$ is the product of the weights of the edges in $M$:
$$
\weight\of M\defeq\prod_{e\in M} \weight\of e.
$$

The \EM{generating function} of  perfect matchings of some graph $G$ with edge weight $\weight$ is
defined as
$$
m\of{G,\weight}\defeq\sum_{M}\weight\of M,
$$
where $M$ ranges over all perfect matchings of $G$.

\ueberschrift{Sign of perfect matchings}
\def\image{\operatorname{img}}
Assuming a fixed proper embedding $\eta$ of graph $G$, the \EM{sign} of some perfect matching
$M$ of $G$ in the embedding $\eta$ is defined as
$$
\sgn\of{M,\eta}\defeq\pas{-1}^{C\of{M,\eta}},
$$
where $C\of{M,\eta}$ is the \EM{number of all crossings of all edges} in $M$,
i.e., denoting by $\image\of{e}$ the image of the (edge--curve corresponding to) edge $e$:
\begin{equation}
\label{eq:crossing}
C\of{M,\eta}\defeq\sum_{\setof{e_1,e_2}\subseteq M} \absof{\image\of{e_1}\cap\image\of{e_2}}.
\end{equation}

It is possible to modify some proper embedding $\eta$ to another embedding $\eta\strich$ so that
the signs of all perfect matchings are the same for $\eta$ and $\eta\strich$: See \figref{fig:crossing}
for a simple example of such \EM{sign--preserving} modification.

\ueberschrift{Signed generating function of perfect matchings}

The \EM{signed generating function} of perfect matchings
$s\of{G,\weight,\eta}$ (which depends not only on the graph $G$, but also on the proper embedding $\eta$)
is defined as
$$
s\of{G,\weight,\eta}\defeq\sum_{M}\sgn\of{M,\eta}\cdot\weight\of M,
$$
where $M$ ranges over all perfect matchings of $G$.

Clearly, $m\of{G,\weight}=s\of{G,\weight,\eta}\equiv 0$ if
$\absof{V\of G}$ is {odd}.

\begin{figure}
\caption{Untangling of a multiple crossing.} 
\label{fig:crossing}
{\small
The left picture shows \EM{four} edge--curves passing through the \EM{same} nonvertex--point: 
Observe that this situation can be ``modified locally'' (in the sense that the embedding is left
unchanged outside some small neighbourhood of the crossing) to a situation like the one
shown in the right picture, and  that such ``local modification''  does not change the sign
of any perfect matching. For instance, if we assume that the edges 
shown here belong to some perfect
matching $M$, both ``local situations'' contribute a factor $\pas{-1}$ to the sign of $M$
for \EM{every two--element subset} of the four curves (i.e., $\pas{-1}^{\binom{4}2}=\pas{-1}^6=1$).
%
}
\begin{center}
\psset{unit=0.5}
\pspicture(-2.5,-3)(11,3)
\psline(-1.5,0)(1.5,0)
\psline[linestyle=dotted](-2.5,0)(-1.5,0)
\psline[linestyle=dotted](1.5,0)(2.5,0)
\psline(1.06066,1.06066)(-1.06066,-1.06066)
\psline[linestyle=dotted](1.06066,1.06066)(1.6,1.6)
\psline[linestyle=dotted](-1.06066,-1.06066)(-1.6,-1.6)
\psline(0.,1.5)(0.,-1.5)
\psline[linestyle=dotted](0,-2.5)(0,-1.5)
\psline[linestyle=dotted](0,1.5)(0,2.5)
\psline(-1.06066,1.06066)(1.06066,-1.06066)
\psline[linestyle=dotted](-1.06066,1.06066)(-1.6,1.6)
\psline[linestyle=dotted](1.06066,-1.06066)(1.6,-1.6)
\drawnonvertex{0,0}
\psline(6,0)(9.5,0)
\psline[linestyle=dotted](5,0)(6,0)
\psline[linestyle=dotted](9.5,0)(10.5,0)
\psline(7.5,2)(7.5,-2)
\psline[linestyle=dotted](7.5,2)(7.5,3)
\psline[linestyle=dotted](7.5,-2)(7.5,-3)
\psline(7,2)(9.5,-0.5)
\psline[linestyle=dotted](7,2)(6.3,2.7)
\psline[linestyle=dotted](9.5,-0.5)(10.2,-1.2)
\psline(9.5,1.5)(6.5,-1.5)
\psline[linestyle=dotted](9.5,1.5)(10.2,2.2)
\psline[linestyle=dotted](6.5,-1.5)(5.8,-2.2)
\drawnonvertex{7.5,0}
\drawnonvertex{8,0}
\drawnonvertex{9,0}
\drawnonvertex{7.5,-0.5}
\drawnonvertex{8.5,0.5}
\drawnonvertex{7.5,1.5}
\endpspicture
\end{center}
\end{figure}

\begin{figure}
\caption{Stembridge's embedding and the Pfaffian.}
\label{fig:pfaffian}
{\small
The left picture shows Stembridge's embedding $\overline{\eta}$ of a graph $G$.

The right picture shows a specific perfect matching $M$ of $G$ for which the number $C\of{M,\overline{\eta}}$
of crossings according to \eqref{eq:crossing} equals $3$, whence the sign
$\sgn\of{M,\overline{\eta}}$ of this matching is $\pas{-1}^3=\pas{-1}$. 
}
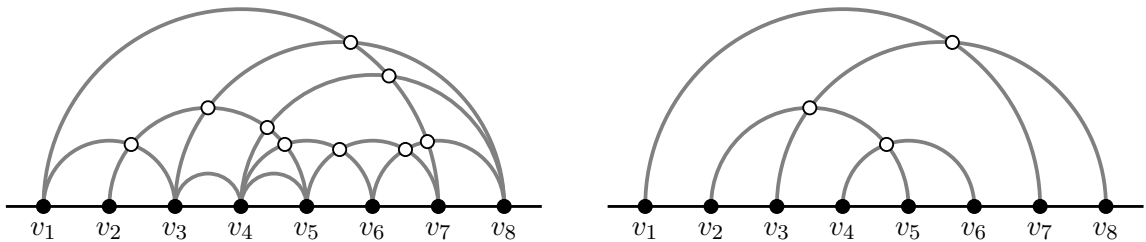
\begin{center}
\input graphics/pfaffians.tex
\end{center}
\end{figure}

\ueberschrift{Stembridge's embedding and the Pfaffian}
For every graph $G$ with vertex set $V\of G=\setof{v_1,v_2,\dots,v_n}$
there is a specific proper embedding which we shall call \EM{Stembridge's embedding}
(see \cite{Stembridge:1990:NPPAPP}):
Every edge $e=\setof{v_i,v_j}$, $i\neq j$, is represented by the half--circle in the upper half--plane
with center $\pas{\frac{i+j}2,0}$ and radius $\frac{\absof{i-j}}2$,
and every vertex $v_i$ is represented by the point $\pas{i,0}$.
%
(See the left picture in \figref{fig:pfaffian} for an example of Stembridge's embedding.)

Denote this specific embedding by $\overline{\eta}$: The \EM{Pfaffian} of a graph $G$ is defined
as the signed generating function for the embedding $\overline{\eta}$, i.e.:
$$
\pfaffian\of{G,\weight}\defeq s\of{G,\weight,\overline{\eta}}.
$$

\ueberschrift{Sign--modifications of edge weights and Kasteleyn's Theorem}
Let $G$ be some simple graph with edge weight $\weight$: Another edge weight $\weight\strich$
is called a \EM{sign--modification} of $\weight$ if for all $e\in E\of G$ there holds
$$
\weight\strich\of e = \weight\of e \text{ or }\weight\strich\of e = -\weight\of e.
$$

Then we may formulate Kasteleyn's Theorem \cite{Kasteleyn:1967:GTACP} as follows:

\begin{thm}[Kasteleyn's Theorem]
\label{thm:kasteleyn}
Let $G$ be a \EM{planar} finite simple graph 
with edge weight $\weight$. Let $\eta$ be an arbitrary proper 
embedding of $G$. Then there exists a sign--modification $\weight\strich$
of $\weight$ such that
$$
m\of{G,\weight} = 
\pfaffian\of{G,\weight\strich}.
$$
\end{thm}

Note that for any \EM{planar} embedding $\eta$ of  $G$ we have $m\of{G,\weight}=s\of{G,\weight,\eta}$.

We will rephrase Speyer's elegant argument \cite{Speyer:2016:VOATOKW} to
give a sketch of proof for the following slight generalization:
\begin{thm}
\label{thm:kasteleyn2}
Let $G$ be a finite simple graph 
with edge weight $\weight$. Let $\eta$ be a 
proper embedding of $G$.
Then for \EM{every} proper embedding $\pas{G,\eta\strich}$ 
there is a sign--modification $\weight\strich$ of $\weight$ such that
$$
s\of{G,\weight,\eta}\equiv s\of{G,\weight\strich,\eta\strich}.
$$
\end{thm}

\section{Speyer's argument, in simple pictures}
The simple idea for the proof of Theorem~\ref{thm:kasteleyn2} is to \EM{successively transform} the
embedding $\eta$
\bit
\item by \EM{local modifications of the embedding}
\item and \EM{corresponding modifications of the edge weight} (if necessary)
\eit
such that the signed generating function of perfect matchings stays \EM{unchanged},
until the embedding $\eta\strich$ is obtained.

Note that smooth bijections $\R^2\to\R^2$ will transform the embedding without changing the (numbers
of) crossings, thus leaving the signed generating function of perfect matchings unchanged.
More informally: If we think of the embedding
as a ``web of infinitely thin and ductile strings'' (corresponding to the edges) glued together
at their endpoints (corresponding to the vertices), then we can imagine that we may ``deform and
drag around'' these strings, and such transformation will not change the signed generating function of
perfect matchings if we do not \EM{remove} crossings or \EM{introduce} new ones.

Moreover, there are modifications which leave the sign of every perfect matching unchanged,
but cannot be achieved by a continuous mapping $\R^2\to\R^2$: One example is the ``untangling of
a multiple crossing''
illustrated in \figref{fig:crossing}.

In the following, we shall consider ``local modifications of proper embeddings''. For simplicity, we assume that
every crossing
is the intersection of \EM{precisely two} edge--curves
(which we can achieve by the ``untangling of
a multiple crossing'' illustrated in \figref{fig:crossing}; without changing the signed generating function
we are interested in).

\ueberschrift{Crossing edges incident to the the same vertex}
If two edges $e_1$, $e_2$ have a vertex $v$ in common, then they can never both belong to
the \EM{same} perfect matching, hence crossings of $e_1$, $e_2$ are
irrelevant for the sign of any perfect matching. If the curve--segments between
$v$ and some crossing $p$ do not contain any other crossing, then $p$
can simply be removed (or introduced)
by the modification illustrated in the following picture:
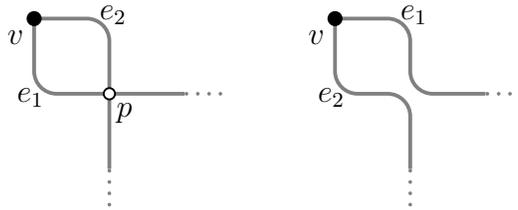
\begin{center}
\psset{unit=0.5}
\pspicture(-2,1)(11,6)
\psline(0,2)(0,6)(-2,6)
\psline(-2,6)(-2,4)(2,4)
\psline[linestyle=dotted](0,1)(0,2)
\psline[linestyle=dotted](2,4)(3,4)
\drawnonvertex{0,4}
\drawvertex{-2,6}
\rput[c](-2.5,5.5){$v$}
\rput[c](-2.1,3.9){$e_1$}
\rput[c](0.1,6.1){$e_2$}
\rput[c](0.4,3.5){$p$}
\psline(8,2)(8,4)(6,4)(6,6)
\psline(6,6)(8,6)(8,4)(10,4)
\psline[linestyle=dotted](8,1)(8,2)
\psline[linestyle=dotted](10,4)(11,4)
\drawvertex{6,6}
\rput[c](5.5,5.5){$v$}
\rput[c](5.9,3.9){$e_2$}
\rput[c](8.1,6.1){$e_1$}
\endpspicture

\end{center}
This modification has no effect on the sign of any perfect matching. If we modify the situation
``from left to right'', we shall call this \EM{straightening out a single crossing} of $e_1$ with $e_2$.

\ueberschrift{Self--crossing edges}
We ruled out self--intersections of edge--curves for proper embeddings, but they might arise by one
of the modifications described below. The following picture makes clear that we can ``change the
situation locally'' and remove (or introduce) such self--crossing.
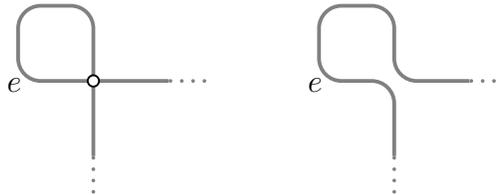
\begin{center}
\psset{unit=0.5}
\pspicture(-2,1)(11,6)
\psline(0,2)(0,6)(-2,6)(-2,4)(2,4)
\psline[linestyle=dotted](0,1)(0,2)
\psline[linestyle=dotted](2,4)(3,4)
\drawnonvertex{0,4}
\rput[c](-2.1,3.9){$e$}
\psline(8,2)(8,4)(6,4)(6,6)(8,6)(8,4)(10,4)
\psline[linestyle=dotted](8,1)(8,2)
\psline[linestyle=dotted](10,4)(11,4)
\rput[c](5.9,3.9){$e$}
\endpspicture

\end{center}
This modification has no effect on the sign of any perfect matching, since only crossings of \EM{different}
edges contribute to the sign.
 If we modify the situation
``from left to right'', we shall call this \EM{straightening out a single self--crossing} of $e$.

\ueberschrift{Dragging a segment of an edge--curve over another segment}
Assume that segments of edges
$e$, $e\strich$ can be ``dragged over one another''
such that precisely two crossings arise (or vanish), as the following picture illustrates:
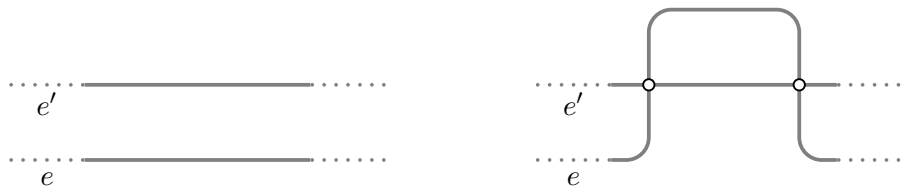
\begin{center}
\psset{unit=0.5}
\pspicture(-12.,-3.)(12.,3)
\psline(-4.,-2.)(-10,-2)
\psline(-4.,0.)(-10,0)
\psline[linestyle=dotted](-4.,-2.)(-2,-2)
\psline[linestyle=dotted](-10,-2)(-12,-2.)
\psline[linestyle=dotted](-4.,0)(-2,0)
\psline[linestyle=dotted](-10,0)(-12,0)
\rput[c](-11, -2.5){{\small $e$}}
\rput[c](-11, -0.5){{\small $e\strich$}}
\psline(4.,-2.)(5.,-2.)(5,2)(9,2)(9,-2)(10,-2)
\psline(4.,0.)(10,0)
\psline[linestyle=dotted](4.,-2.)(2,-2)
\psline[linestyle=dotted](10,-2)(12,-2.)
\psline[linestyle=dotted](4.,0)(2,0)
\psline[linestyle=dotted](10,0)(12,0)
\drawnonvertex{5.,0}
\drawnonvertex{9.,0}
\rput[c](3, -2.5){{\small $e$}}
\rput[c](3, -0.5){{\small $e\strich$}}
\endpspicture

\end{center}
This modification has  no effect on the sign of \EM{any} perfect matching.
 If we modify the situation
``from left to right'', we shall call this \EM{straightening out a double crossing} of $e$ with $e\strich$.

\ueberschrift{Transition of a vertex through an edge}
Assume that a segment of 
some edge $e$ can be ``dragged over some vertex $v$''
such that precisely one crossing arises (or vanishes) for every edge incident with $v$,
as the following picture illustrates:
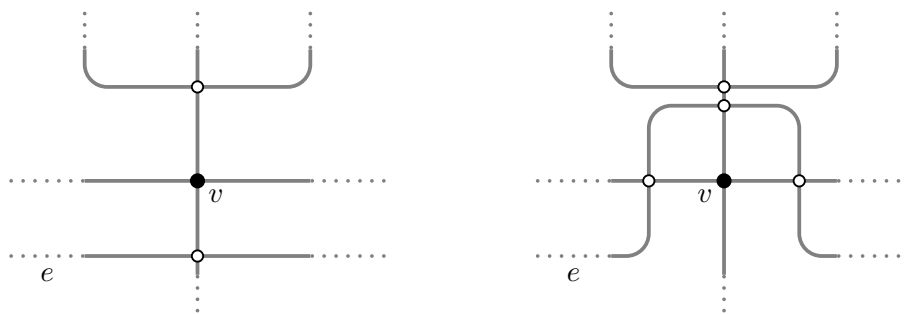
\begin{center}
\psset{unit=0.5}
\pspicture(-12.,-3.5)(12.,5)
\psline(-4.,-2.)(-10,-2)
\psline(-4.,0.)(-10,0)
\psline(-7,-2.5)(-7,3.5)
\psline[linestyle=dotted](-4.,-2.)(-2,-2)
\psline[linestyle=dotted](-10,-2)(-12,-2.)
\psline[linestyle=dotted](-4.,0)(-2,0)
\psline[linestyle=dotted](-10,0)(-12,0)
\psline[linestyle=dotted](-7,3.5)(-7.,4.5)
\psline[linestyle=dotted](-7,-2.5)(-7.,-3.5)
\psline(-4,3.5)(-4,2.5)(-10,2.5)(-10,3.5)
\psline[linestyle=dotted](-4.,3.5)(-4,4.5)
\psline[linestyle=dotted](-10.,3.5)(-10,4.5)
\drawvertex{-7,0}
\drawnonvertex{-7.,2.5}
\drawnonvertex{-7.,-2}
\rput[c](-11, -2.5){{\small $e$}}
\rput[c](-6.5, -0.4){{\small $v$}}
\psline(4.,-2.)(5.,-2.)(5,2)(9,2)(9,-2)(10,-2)
\psline(4.,0.)(10,0)
\psline(7,-2.5)(7,3.5)
\psline(4,3.5)(4,2.5)(10,2.5)(10,3.5)
\psline[linestyle=dotted](4.,3.5)(4,4.5)
\psline[linestyle=dotted](10.,3.5)(10,4.5)
\psline[linestyle=dotted](4.,-2.)(2,-2)
\psline[linestyle=dotted](10,-2)(12,-2.)
\psline[linestyle=dotted](4.,0)(2,0)
\psline[linestyle=dotted](10,0)(12,0)
\psline[linestyle=dotted](7,3.5)(7.,4.5)
\psline[linestyle=dotted](7,-2.5)(7.,-3.5)
\drawvertex{7,0}
\drawnonvertex{5.,0}
\drawnonvertex{7.,2}
\drawnonvertex{9.,0}
\drawnonvertex{7.,2.5}
\rput[c](3, -2.5){{\small $e$}}
\rput[c](6.5, -0.4){{\small $v$}}
\endpspicture

\end{center}
Since every perfect matching $M$ must contain \EM{precisely one} edge incident to $v$, the sign of $M$
is \EM{reversed} by this modification if and only if $e\in M$: But we can easily offset this by replacing
$\weight\of e$ by $-\weight\of e$.

\section{Sketch of proof by example}
\label{sec:example}

A rigorous proof of Theorem~\ref{thm:kasteleyn2} 
would involve tedious technical details: We shall avoid them by 
only illustrating the \EM{idea of proof}, i.e, the ``successive modification of embeddings'', 
in an example.

The example we shall consider is the \EM{complete bipartite graph} $K_{3,3}$: Three
different proper embeddings of this graph are shown in \figref{fig:proper}.

Observe that the embedding shown in the middle of \figref{fig:proper} can easily
be transformed to the embedding
shown in the right by the ``untangling of a multiple crossing''
illustrated in \figref{fig:crossing}.

\figref{fig:example} illustrates the process of successively modifying the embedding $\eta$
shown in the left picture  of \figref{fig:proper} until the embedding $\eta\strich$ shown in the middle
picture of \figref{fig:proper} is obtained.
We shall explain the single steps in the following, where we denote the edge connecting
vertices $v_i$ and $v_j$ by $e_{i,j}$.

\begin{figure}
\caption{Three different embeddings of the same graph $G$.}
\label{fig:proper}
{\small
The graph $G$ we are considering here is the \EM{complete bipartite graph} $K_{3,3}$:
The fact that $K_{3,3}$ is \EM{not planar} corresponds to the fact that we are
unable to present an embedding without crossings.
}
\begin{center}
\input graphics/proper.tex
\end{center}
\end{figure}

\begin{figure}
\caption{Example: Succesive modification of embeddings.}
\label{fig:example}
{\small
Consider the embeddings of the graph $G$ depicted in \protect{\figref{fig:proper}}:
We want to transform the embedding $\eta$ shown in the left picture to the embedding $\eta\strich$
shown in the middle picture of \protect{\figref{fig:proper}}.
}
\begin{center}
\input graphics/example.tex
\end{center}
\end{figure}

The left upper picture in \figref{fig:example} shows the result of a smooth deformation of the plane,
where vertices $v_2$, $v_4$ and $v_6$ already have arrived at their desired positions
(according to the hexagonal configuration in $\eta\strich$). Now we want to drag $v_5$ to its 
desired position along the path indicated by the dashed arrow. This involves a transition
of $v_5$ through edge $e_{3,6}$, whose weight must therefore change its sign.

The upper picture in the middle shows the result of this operation: Note that we may now
straighten out the double crossing of $e_{4,5}$ with $e_{3,6}$ and the single
crossing of $e_{5,6}$ with  $e_{3,6}$.

Now we want to drag $v_1$ to its 
desired position along the path indicated by the dashed arrow. This involves a transition
of $v_1$ through edge $e_{2,3}$, whose weight must therefore change its sign. Note that after performing
this operation, we may straighten out the double crossing of 
$e_{1,4}$ with $e_{2,3}$.

The right upper picture shows the result of these operations: Observe that we may now straighten
out the double crossing of $e_{2,3}$ with $e_{1,6}$ and the single crossing of $e_{1,2}$ with $e_{2,3}$.

Now we want to drag $v_3$ to its 
desired position along the path indicated by the dashed arrow: This can be done without introducing
or removing any crossing.

The left lower picture shows the result of these operations. By dragging $e_{1,6}$ over vertex
$v_4$ and straightening out the single crossing with $e_{1,4}$ thus introduced, we arrive at the lower
picture in the middle: This involves a transition
of $v_4$ through $e_{1,6}$, whose weight must therefore change its sign. 

Now dragging $e_{3,6}$ over vertices $v_4$ \EM{and} $v_5$ (and changing the sign of $\weight\of{e_{3,6}}$
\EM{twice}, accordingly), and straightening out
\bit
\item the two single crossings of $e_{3,6}$ with $e_{5,6}$ and $e_{3,4}$
\item and the double crossing of $e_{3,6}$ with $e_{4,5}$
\eit
introduced by this modification gives the right lower picture. 

Repeating the last step by dragging $e_{2,5}$ over vertices $v_3$ and $v_4$, we arrive at the right proper drawing
in {\figref{fig:proper}.

\bibliography{/Users/mfulmek/Work/Forschung/Basis/database}

\end{document}